\documentclass[12pt]{article}
\usepackage{amscd}
\usepackage{mathrsfs}
\usepackage{amsfonts}
\usepackage{graphicx}
\usepackage{amsmath,amscd,stmaryrd,amssymb,array, amsfonts,amsmath,amssymb}
\usepackage{enumitem}
\usepackage{amsthm}
\usepackage{bm}

\setenumerate[1]{itemsep=0pt,partopsep=0pt,parsep=\parskip,topsep=5pt}
\setitemize[1]{itemsep=0pt,partopsep=0pt,parsep=\parskip,topsep=0pt}
\setdescription{itemsep=0pt,partopsep=0pt,parsep=\parskip,topsep=5pt}

\numberwithin{figure}{section}
 \numberwithin{equation}{section}

\newtheorem{theorem}{Theorem}[section]
\newtheorem{proposition}[theorem]{Proposition}
\newtheorem{definition}[theorem]{Definition}

\newtheorem{lemma}[theorem]{Lemma}
\newtheorem{remark}[theorem]{Remark}

\newcommand{\bR}{{\bf R}}

\def\be{\begin{equation}}
\def\ee{\end{equation}}
\def\bes{\begin{equation*}}
\def\ees{\end{equation*}}
\def\bsp{\begin{split}}
\def\esp{\end{split}}

\def\ba{\begin{array}}
\def\ea{\end{array}}
\def\benu{\begin{enumerate}}
\def\eenu{\end{enumerate}}
\def\bt{\begin{theorem}}
\def\et{\end{theorem}}
\def\bp{\begin{proposition}}
\def\ep{\end{proposition}}
\def\bl{\begin{lemma}}
\def\el{\end{lemma}}
\def\br{\begin{remark}}
\def\er{\end{remark}}
\def\bd{\begin{definition}}
\def\ed{\end{definition}}

\def\W{\Omega}

\def\.{\cdot}

\def\~{\tilde}
\def\8{\infty}

\def\Vs{\vskip8pt}\def\vs{\vskip4pt}

\def\({\left(}\def\){\right)}

\begin{document}

\begin{center}
{\bf\Large Notes on A Superlinear Elliptic Problem}
\end{center}

\vs\centerline{Haoyu  Li
%\footnote{Supported by the grant of NSFC (10771159,
%11071185)}
}  %\vskip20pt
\begin{center}
{\footnotesize
{Center of Applied Mathematics,  Tianjin University\\
          Tianjin 300072,  China\\

{\em E-mail}:  hyli1994@hotmail.com
}}
\end{center}

\Vs

{\footnotesize
{\bf Abstract.}  Based on the theory of invariant sets of descending flow, we give a new proof of the existence of three nontrivial solutions and some remarks on it.
 \Vs
\Vs {\bf 2010 MSC:} 34C23, 34K18, 35B32, 37G99.

%\Vs {\bf Running Head:} Dynamic Bifurcation of Evolution Equations.

%\Vs\Vs {\bf Date:}   June 16-{\em th}, 2015.

}

\tableofcontents

\section{Introduction}

The existence and multiplicity of solutions for superlinear elliptic problem are classical problems in the field of variational methods. It is well-known that under some general assumptions, we can obtain a positive and a negative solution. If in addition, let the nonlinear term be odd, we will have infinitely many solutions. These are classical results which can be found in \cite{AR} and \cite{R}. In 1991, Wang discovered a third solution in \cite{W} via the local linking and Morse theory. In 2001, Liu and Sun proved the existence of four solutions in \cite{LS} under a general framework of invariant sets of descending flow. One positive solution, one negative solution, one sign-changing solution and one possible trivial solution. The effort on refining this result is never ceased. A classical improved version is \cite{LW}, which can obtain more solutions with prescribed sign with some additional assumptions. In a recent work, \cite{LL}, a result for four nontrivial solutions is given with a finer analysis.

In this note we give a new proof of the classical theorem on the existence of three solutions and some remarks on it. Let us make the following assumptions:
\begin{itemize}
  \item [$(f1)$] The function $f\in C(\mathbb{R})$ satisfies
  \begin{align}
  f(u)&=o(|u|),\,\,as\,\,u\to0;\nonumber\\
  f(u)&=O(|u|^{p-1}),\,\,as\,\,u\to\infty,\nonumber
  \end{align}
  where $p\in(2,2^{*})$;
  \item [$(f2)$] There is a number $\mu>2$ with $0<\mu F(u)\leq uf(u)$, where
  $$F(u)=\int_{0}^{u}f(s)ds.$$
\end{itemize}
Consider the superlinear elliptic problem
\begin{equation}\label{e:A}
\left\{
\begin{aligned}
-\Delta u=f(u),\quad in \,\,\Omega,\\
u=0,\qquad on \,\,\partial\Omega.
\end{aligned}
\right.
\end{equation}
where $\Omega\subset\mathbb{R}^{N}$ is a domain whose boundary regular enough for the Sobolev embeddings and $N\geq3$.
\begin{theorem}\label{t:A}
Suppose that the assumptions $(f1)$ and $(f2)$ are satisfied, the equation (\ref{e:A}) possesses at least three solutions, a positive solution, a negative solution and a sign-changing solution.
\end{theorem}
\br
 Although this is a classical theorem with a large number of proofs, to the best of the author's knowledge, the proof we give in this note is new. We refer the \cite{LLW}, \cite{LS}, and \cite{W} for a historical reference. For more results on superlinear elliptic problems, readers can find them in \cite{BWW}.
\er

The structure of this paper is as follows: In the second section, we will give some notations, definitions and lemmas. And we will prove the existence of critical points in the third section. In the last section, we give some remarks on the more constructions of the sign-changing minimax values and the linking sets.

\section{Preliminaries}
\subsection{Notations}
\begin{itemize}
  \item An open $\delta-$neighbourhood of a set $A$ is denoted by $A_{\delta}$;
  \item The complementary set of the set $A$ is denoted by $A^{c}$;
  \item For a subset $A$ of a linear space, let $-A:=\{-x|x\in A\}$;
  \item $E:= H_{0}^{1}(\Omega)$ denotes the usual Sobolev space. By $\|\cdot\|$, we denote the Sobolev norm and $|\cdot|_{p}$ the $L^{p}$ norm. The closure in this note is always taken with respect to the Sobolev norm;
  \item We denote the open ball in the Sobolev space $H_{0}^{1}(\W)$ with radii $\rho$ and center at $0\in H_{0}^{1}(\W)$ by $B_{\rho}$, by $\partial B_{\rho}$ its boundary;
  \item $P^{\pm}=\{u\in H_{0}^{1}(\Omega)|u\geq (\leq) 0\,\,a.e.\,\,in\,\,\Omega\}$ is positive (negative) cone in the Sobolev space $H_{0}^{1}(\W)$;
  \item $W_{\varepsilon}:=\overline{P_{\varepsilon}^{+}}\cup\overline{P_{\varepsilon}^{-}}$;
  \item $I^{a}=\{u\in E|I(u)\leq a\}$ and $I^{a}_{b}:=I^{a}\backslash I^{b}$ are level sets;
  \item $K=\{u\in E|I'(u)=0\}$, $K_{a}=\{u\in E|I'(u)=0\,\,and\,\,I(u)=a\}$ and $K^{a}_{b}:=K\cap I^{a}_{b}$ are the sets of critical points of $I$ in $E$.
\end{itemize}
\subsection{Definitions and lemmas}
The following definitions and lemmas are standard for modern variational methods for sign-changing critical points. We still state it here for the completeness. We refer \cite{LLW}, \cite{LS} and \cite{LWZ} for a general theoretic construction.

\begin{definition}\label{d:invariant}
Let $X$ be a Banach space and $I\in C^{1}(X,\mathbb{R})$. Then $P\subset X$ is called an admissible invariant set with respect to $I$ at level $c$, if $K_{c}\backslash P=\emptyset$ implies the existence of a positive constant $\varepsilon_{0}$ and a function $\eta\in C(E,E)$ with
\begin{itemize}
  \item[$(1)$] $\eta(\overline{P})\subset\overline{P}$,
  \item[$(2)$] $\eta|_{I^{c-\varepsilon}}=id_{I^{c-\varepsilon}}$,
  \item[$(3)$] $\eta(I^{c+\varepsilon} \backslash P)\subset I^{c-\varepsilon}$,
\end{itemize}
for any $\varepsilon\in(0,\varepsilon_{0})$.
\end{definition}

\br
There is a symmetric version of definition \ref{d:invariant}. Readers can find it in \cite{LLW} and \cite{LWZ}.
\er

The energy of the Problem (\ref{e:A}) is
$$I(u)=\frac{1}{2}\int_{\W}|\nabla u|^{2}-\int_{\W} F(u),$$
for $u\in H_{0}^{1}(\Omega)=E$. It should be noted that the $I$ satisfies the $(PS)$ condition. This is a classical result. We will claim Lemma \ref{l:PS} and refer \cite{W} for its proof.
\begin{lemma}\label{l:PS}
$I$ satisfies the $(PS)$ condition.
\end{lemma}
In the following part, we shall state some basic properties of $\overline{P_{\varepsilon}^{+}}$, $\overline{P_{\varepsilon}^{-}}$ and $W_{\varepsilon}$. Consider the operator $A$ defined by
$$v=(-\Delta)^{-1}(f(u))=:A(u).$$
Then the operator $A:E\to E$ is well-defined, continuous and compact. One can find an analogue in \cite[Section 3.1]{LWZ}.

\begin{lemma}
\begin{itemize}\label{l:A}
  \item[$(1)$] $I'(u)(u-A(u))=\|u-A(u)\|^{2}$,
  \item[$(2)$] $\|I'(u)\|=\|u-A(u)\|.$
\end{itemize}
\end{lemma}
\noindent{\bf Proof.}
\begin{itemize}
  \item[$(1)$] Since $u,\,A(u)\in H_{0}^{1}(\W)$,
    \begin{align}
    I'(u)(u-A(u))&=\int_{\W} (-\Delta)(u-A(u))(u-A(u))\nonumber\\
                          &=\int_{\W} |\nabla(u-A(u))|^{2}\nonumber\\
                          &=\|u-A(u)\|^{2}.\nonumber
    \end{align}
  \item[$(2)$]For any $\varphi\in H_{0}^{1}(\Omega)$,
    \begin{align}
     I'(u)\varphi &=(u-A(u),\varphi)_{H}\nonumber\\
                           &\leq\|u-A(u)\|\cdot\|\varphi\|.\nonumber
    \end{align}
\end{itemize}
Hence $\|I'(u)\|\leq\|u-A(u)\|$. Similarly,
$$\|I'(u)\|\geq\|u-A(u)\|.$$
\begin{flushright}
\qedsymbol
\end{flushright}

Since every solution can be considered as a fixed point of the mapping (and hence of the flow), we have the following lemma.

\begin{lemma}\label{l:attract}
There is an $\varepsilon_{0}>0$ such that for any $\varepsilon\in(0,\varepsilon_{0})$, $A(\overline{P_{\varepsilon}^{\pm}})\subset \overline{P_{\varepsilon}^{\pm}}$ and every nontrivial solution in $\overline{P_{\varepsilon}^{\pm}}$ is positive (negative) i.e. the solution belongs to $P_{\varepsilon}^{\pm}$.
\end{lemma}
This lemma is similar to \cite[Lemma 3.4]{LWZ}. For the completeness, we give the proof here.

\noindent{\bf Proof.}
We only proof the theorem for $P^{-}_{\varepsilon}$. By $(g1)$ and $(g2)$, for any $\delta>0$, there exists a constant $C_{\delta}>0$ such that
$$|f(u)|\leq\delta|u|+C_{\delta}|u|^{q-1}.$$
Let $u\in E$ and $v=A(u)$. Since for any $r\in[2,2^{*}]$, there is a constant $m_{r}>0$ such that
$$|u^{\pm}|_{r}=\inf_{w\in P^{\mp}}|u-w|_{r}\leq m_{r}\inf_{w\in P^{\mp}}\|u-w\|=m_{r}d(u,P^{\mp}).$$
Combining the fact that $d(v,P^{-})\leq\|v^{+}\|$ and $(g2)$,
\begin{align}
d(v,P^{-})\|v^{+}\|&\leq\|v^{+}\|^{2}=(v,v^{+})_{H}\nonumber\\
                   &=\int_{\W}f(u)v^{+}\leq\delta\int_{\W}|u^{+}|+C_{\delta}\int_{\W}|u^{+}|^{p-1}|v^{+}|\nonumber\\
                   &\leq\delta|u^{+}|_{2}|v^{+}|_{2}+C_{\delta}|u^{+}|_{p}^{p-1}|v^{+}|_{p}\nonumber\\
                   &\leq C(\delta d(u,P^{-})+C_{\delta}d(u,P^{-})^{p-1})\|v^{+}\|.\nonumber
\end{align}
It follows that
$$d(A(u),P^{-})\leq C(\delta d(u,P^{-})+C_{\delta}d(u,P^{-})^{p-1}).$$
Choosing $\delta$ small enough, there exists a $\varepsilon_{0}>0$ such $\varepsilon\in(0,\varepsilon_{0})$,
$$d(A(u),P^{-})\leq\frac{1}{2}d(u,P^{-})$$
for any $u\in \overline{P_{\varepsilon}^{-}}$, which implies that $A(\overline{P_{\varepsilon}^{-}})\subset\overline{P_{\varepsilon}^{-}}$. Then if there is any point $u\in \overline{P_{\varepsilon}^{-}}$ such that $A(u)=u$, we will have $u\in P^{-}$. Using the maximum principle, if $u\neq\theta$, then $u<0$ a.e. in $\W$.
\begin{flushright}
\qedsymbol
\end{flushright}

\begin{lemma}\label{l:B}
There is a Lipschitz mapping $B:E\backslash K_{c}\to E$ such that
\begin{itemize}
  \item[$(1)$] $B(\overline{P_{\varepsilon}^{\pm}})\subset \overline{P_{\varepsilon}^{\pm}}$ for\,\,$0<\varepsilon<\varepsilon_{0}$;
  \item[$(2)$] $\frac{1}{2}\|u-B(u)\|\leq\|u-A(u)\|\leq2\|u-B(u)\|$ for any $u\in H_{0}^{1}(\Omega)$;
  \item[$(3)$] $I'(u)(u-B(u))\geq\frac{1}{2}\|u-A(u)\|^{2}$ for any $u\in E\backslash K$;
\end{itemize}
\end{lemma}

This is \cite[Lemma 2.1]{BLW}. The main idea of this lemma is to replace the gradient field by its pseudo-gradient field. We refer \cite[Lemma 2.1]{BLW} for the proof.

\begin{lemma}\label{l:mapping}
If for some $\varepsilon>0$ small enough, $K_{c}\backslash W_{\varepsilon}=\phi$, then there is a $\varepsilon_{0}>0$, for $0<\varepsilon<\varepsilon '<\varepsilon_{0}$, there is a continuous mapping $\eta:[0,1]\times H_{0}^{1}(\Omega)\to H_{0}^{1}(\Omega)$ such that
\begin{itemize}
  \item[$(1)$] $\eta(0,u)=u$ for any $u\in H_{0}^{1}(\Omega)$.
  \item[$(2)$] $\eta(t,u)=u$ for all $t\in[0,1]$,$\,\,u\notin I^{c+\varepsilon '}_{c-\varepsilon '}$.
  \item[$(3)$] $\eta(1,I^{c+\varepsilon}\backslash W_{\varepsilon})\subset I^{c-\varepsilon}$.
  \item[$(4)$] $\eta(t,\overline{P_{\varepsilon}^{\pm}})\subset\overline{P_{\varepsilon}^{\pm}}$ for all $t\in[0,1]$.
\end{itemize}
\end{lemma}
\noindent{\bf Proof.}
This proof is similar to the one of \cite[Lemma 3.6]{LWZ}. But we give it here for the completeness. We denote $W_{\varepsilon}$ by $W$ for short.

Since $K_{c}\backslash W=\emptyset$, there is a $\delta>0$ such that $(K_{c})_{\delta}\subset W$. Since $I$ satisfies the (PS)-condition, there exist $\varepsilon_{0},\alpha>0$ such that
$$\|I'(u)\|\geq\alpha,$$
for $u\in I^{-1}([c-\varepsilon_{0},c+\varepsilon_{0}])\backslash(K_{c})_{\frac{\delta}{2}}$. By Lemma \ref{l:A} and Lemma \ref{l:B}, there is a constant $\beta>0$ such that
$$\|u-B(u)\|\geq\alpha,$$
for $u\in I^{-1}([c-\varepsilon_{0},c+\varepsilon_{0}])\backslash(K_{c})_{\frac{\delta}{2}}$. Without loss of generality, assume that $\varepsilon_{0}\leq\frac{\beta\delta}{32}$. Define
$$V(u):=\frac{u-B(u)}{\|u-B(u)\|}$$
for $u\in E_{0}=E\backslash K$ and a Lipschitz function $g:E\to[0,1]$ such that
\begin{equation}\label{e:cutoff}
g(u)=\left\{
\begin{aligned}
0,\,\, & if\,\, u\notin I_{c-\varepsilon'}^{c+\varepsilon'}\,\, or \,\,u\in(K_{c})_{\frac{\delta}{4}}, \\
1,\,\, & if\,\, u\in I_{c-\varepsilon'}^{c+\varepsilon'}\,\, or \,\,u\notin(K_{c})_{\frac{\delta}{2}}.
\end{aligned}
\right.
\end{equation}
Consider the Cauchy problem
\begin{equation}\label{e:flow}
\left\{
\begin{aligned}
&\frac{d\tau}{dt}=-g(\tau)V(\tau),\\
&\tau(0,u)=u.
\end{aligned}
\right.
\end{equation}
For any $u\in E$, the problem \ref{e:flow} admits a unique solution $\tau(\cdot,u)\in C(\bR^{+},E)$. Define $\eta(t,u)=\tau(\frac{16\varepsilon}{\beta}t,u)$. And hence (1) and (2) are proved.

Let $u\in I^{c+\varepsilon}\backslash W$. $I(\tau(t,u))$ is decreasing for $t\geq0$. If there is a $t_{0}\in[0,\frac{16\varepsilon}{\beta}]$ with $I(\tau(t_{0},u))<c-\varepsilon$ then $I(\tau(1,u))<c-\varepsilon$. Otherwise, for any $t\in[0,\frac{16\varepsilon}{\beta}]$, $I(\tau(t,u))\geq c-\varepsilon$. Then, $\tau(t,u)\in I^{-1}[c-\varepsilon,c+\varepsilon]$ for $t\in[0,\frac{16\varepsilon}{\beta}]$. We claim that for any $t\in[0,\frac{16\varepsilon}{\beta}]$, $\tau(t,u)\notin(K_{c})_{\frac{\delta}{2}}$. If $\tau(t_{1},u)\in(K_{c})_{\frac{\delta}{2}}$
 for some $t_{1}\in[0,\frac{16\varepsilon}{\beta}]$, then we have
$$\frac{\delta}{2}\leq\|\tau(t_{1},u)\|\leq\int_{0}^{t_{1}}\|\do{\tau}(s,u)\|ds\leq t_{1}\leq\frac{16\varepsilon}{\beta}$$
since $u\notin(K_{c}))_{\frac{\delta}{2}}$. And we have a contradiction with $\varepsilon<\varepsilon_{0}\leq\frac{\beta\delta}{32}$. So $g(\tau(t,u))=1$ for $t\in[0,\frac{16\varepsilon}{\beta}]$. Then, by Lemma \ref{l:B},
\begin{align}
I\Big(\tau\big(\frac{16\varepsilon}{\beta},u\big)\Big) &= I(u)-\int_{0}^{\frac{16\varepsilon}{\beta}}\big<I'(\tau(s,u)),V(\tau(s,u))\big>ds\nonumber\\
& \leq I(u)-\frac{1}{8}\int_{0}^{\frac{16\varepsilon}{\beta}}\|\tau(s,u)-B\tau(s,u)\|ds\nonumber\\
& \leq c-\varepsilon.\nonumber
\end{align}

And one can find the proof of (4) in \cite{LS}.
\begin{flushright}
\qedsymbol
\end{flushright}
Let $\overline{\eta}(u)=\eta(1,u)$. And hence $P_{\varepsilon}^{\pm}$ are admissible invariant set with respect to $I_{\theta}$ at level $c$ when $K_{c}\backslash P_{\varepsilon}^{\pm}=\emptyset$. The same with $W$.

Following the same idea, by modifying the equation (\ref{e:cutoff}), we can proof the existence of the next mapping.
\begin{lemma}\label{l:mapping2}
For some $0<\varepsilon_{1}<\varepsilon_{2}$, if $K_{c}\backslash W_{\varepsilon_{2}}=\emptyset$, then there is a $\varepsilon_{0}>0$, for $0<\varepsilon<\varepsilon '<\varepsilon_{0}$, there is a continuous mapping $\eta:[0,1]\times H_{0}^{1}(\Omega)\to H_{0}^{1}(\Omega)$ such that
\begin{itemize}
  \item[$(1)$] $\eta(0,u)=u$ for any $u\in H_{0}^{1}(\Omega)$.
  \item[$(2)$] $\eta(t,u)=u$ for all $t\in[0,1]$,$\,\,u\in \big(I_{c-\varepsilon '}^{c+\varepsilon '}\big)^{c}\cup W_{\varepsilon_{1}}$.
  \item[$(3)$] $\eta(1,I^{c+\varepsilon}\backslash W_{\varepsilon_{2}})\subset I^{c-\varepsilon}$.
  \item[$(4)$] $\eta(t,\overline{P_{\varepsilon_{2}}^{\pm}})\subset\overline{P_{\varepsilon_{2}}^{\pm}}$ for all $t\in[0,1]$.
\end{itemize}
\end{lemma}

\section{Proof of Theorem \ref{t:A}}
Consider the following three minimax values
$$c_{+}=\inf_{g\in\Gamma_{+}}\sup_{t\in[0,1]}I(g(tRe_{1})),$$
$$c_{-}=\inf_{g\in\Gamma_{-}}\sup_{t\in[0,1]}I(g(-tRe_{1}))$$
and
$$c_{s}=\inf_{g\in\Gamma_{s}}\sup_{u\in g(B_{2}^{+})\backslash W_{\varepsilon_{2}}}I(u),$$
where for the value $c_{+}$, we denote
$$\Gamma_{+}=\{g\in C(E,E)|\,\,g(0)=0,\,g(Re_{1})=Re_{1}\,\,and\,\,g(\overline{P^{+}_{\varepsilon}})\subset \overline{P^{+}_{\varepsilon}}\},$$
$e_{1}$ is the principle eigenfunction with positive sign and unit Sobolev norm and $R>0$ is such a number with $I(Re)<0$ and $I(-Re)<0$; for the value $c_{-}$, we denote
$$\Gamma_{-}=\big\{g\in C(E,E)\big|\,\,g(0)=0,\,g(-Re_{1})=-Re_{1}\,\,and\,\,g(\overline{P^{-}_{\varepsilon}})\subset \overline{P^{-}_{\varepsilon}}\big\},$$
and for the value $c_{s}$, we denote
\begin{itemize}
  \item $\Gamma_{s}=\{g\in C(E,E)|\,\,g|_{\partial B^{+}_{2}}=id\};$
  \item $B_{2}^{+}=\overline{B_{R}}\cap E_{2}^{+};$
  \item $E_{2}^{+}=[0,+\infty)e_{2}\oplus\mathbb{R}e_{1},$
  where $e_{2}$ is a sign-changing function, say, the second eigenvector of $-\Delta$ with zero Dirichlet boundary value;
  \item $\partial_{0}B_{2}^{+}=\partial B_{R}\cap E_{2}^{+}$, $\partial B^{+}_{2}=\partial_{0}B_{2}^{+}\cup \big(\overline{B_{R}}\cap\mathbb{R}e_{1}\big)$;
  \item $W_{\varepsilon_{2}}=\overline{P^{+}_{\varepsilon_{2}}}\cup\overline{P^{-}_{\varepsilon_{2}}}$ for some positive $\varepsilon_{2}$ small enough.
\end{itemize}
Let $R>0$ be a number large enough such that $I|_{\partial_{0}B^{+}_{2}}<0$. We shall prove that
\begin{align}
K_{c_{+}}\cap P^{+}\neq\emptyset,\label{c:A}\\
K_{c_{+}}\cap P^{-}\neq\emptyset,\label{c:B}
\end{align}
and
\begin{align}
K_{c_{s}}\backslash W_{\varepsilon_{2}}\neq\emptyset,\label{c:C}
\end{align}
for some suitable positive $\varepsilon_{2}$.
\subsection{Verification of (\ref{c:A})}
Since (\ref{c:A}) and (\ref{c:B}) are similar, we only prove (\ref{c:A}) here. Due to the Lemma \ref{l:attract}, we only need to prove that
$$K_{c_{+}}\cap \overline{P_{\varepsilon}^{+}}\neq\emptyset$$
for $\varepsilon>0$ small enough. Let us assume that $K_{c_{+}}\cap \overline{P_{\varepsilon}^{+}}=\emptyset$ holds. Using the (PS) condition, we can claim that there is a small positive number $\delta$ such that
$$K^{c_{+}+2\delta}_{c_{+}-2\delta}\cap\overline{P_{\varepsilon}^{+}}=\emptyset.$$
In the following part, we only need to verify that there is a descending flow $\eta$ as in Lemma \ref{l:mapping} satisfies that
$$\eta(I^{c_{+}+\delta}\cap\overline{P_{\varepsilon}^{+}})\subset I^{c_{+}-\delta}\cap\overline{P_{\varepsilon}^{+}}.$$
And the rest part of the proof is a routine, we refer \cite{R} for details.

On one hand, since
$$I^{c_{+}+\delta}\cap\overline{P_{\varepsilon}^{+}}\subset\overline{P_{\varepsilon}^{+}},$$
we have
$$\eta(I^{c_{+}+\delta}\cap\overline{P_{\varepsilon}^{+}})\subset\overline{P_{\varepsilon}^{+}}.$$
On the other hand, we claim that the norm of the pseudo-gradient vector field has a positive lower bound on the closed set $\overline{I^{c_{+}+\delta}_{c_{+}-\delta}}\cap\overline{P_{\varepsilon}^{+}}$. Otherwise, there is a sequence $(u_{n})_{n}\subset\overline{I^{c_{+}+\delta}_{c_{+}-\delta}}\cap\overline{P_{\varepsilon}^{+}}$ such that
$$I'(u_{n})_{n}\to0$$
and
$$(I(u_{n}))_{n}\subset[c_{+}-\delta,c_{+}+\delta]$$
as $n\to\infty$. Using the (PS) condition, there is a $u^{*}\in\overline{I^{c_{+}+\delta}_{c_{+}-\delta}}\cap\overline{P_{\varepsilon}^{+}}$ with $I'(u^{*})=0$. Hence we have a contradiction with $K^{c_{+}+2\delta}_{c_{+}-2\delta}\cap\overline{P_{\varepsilon}^{+}}=\emptyset$. Therefore the homeomorphism $\eta$ defined by the pseudo-gradient vector field satisfies
$$\eta(\overline{I^{c_{+}+\delta}_{c_{+}-\delta}}\cap\overline{P_{\varepsilon}^{+}})\subset I^{c_{+}-\delta},$$
which implies that
$$\eta(I^{c_{+}+\delta}\cap\overline{P_{\varepsilon}^{+}})\subset I^{c_{+}-\delta}\cap\overline{P_{\varepsilon}^{+}}.$$

\subsection{Verification of (\ref{c:C})}
\subparagraph{Step 1}

Firstly, we verify a linking-type result, i.e. for any $g\in \Gamma_{s}$
$$(g(B^{+}_{2})\cap\partial B_{\rho})\backslash W_{\varepsilon_{2}}\neq \emptyset.$$
For any $g\in\Gamma_{s}$, there is a odd mapping $\overline{g}:B_{2}\to E$ defined as
\begin{equation}
\overline{g}(u)=\left\{
\begin{aligned}
g(u),\qquad  & u\in B_{2}^{+}; \nonumber\\
-g(-u),\qquad & -u\in B_{2}^{+}, \nonumber
\end{aligned}
\right.
\end{equation}
where $B_{2}=B_{2}^{+}\cup (-B_{2}^{+})$. Using a genus argument (c.f. \cite[Lemma 6.4]{S1}),
$$\gamma(\partial B_{\rho}\cap\overline{g}(B_{2}))=2$$
and
$$\gamma(\partial B_{\rho}\cap W_{\varepsilon_{2}})=1,$$
for $\varepsilon_{2}$ small enough. So $(\overline{g}(B_{2})\cap\partial B_{\rho})\backslash W_{\varepsilon_{2}}= \emptyset$, i.e. $\partial B_{\rho}\cap\overline{g}(B_{2})\subset\partial B_{\rho}\cap W_{\varepsilon_{2}}$, is not possible. Since $\overline{g}$ is odd, we have
$$(g(B^{+}_{2})\cap\partial B_{\rho})\backslash W_{\varepsilon_{2}}\neq \emptyset.$$

\subparagraph{Step 2}

Secondly, we prove that $K_{c_s}\backslash W_{\varepsilon_{1}}\neq\emptyset$. Since for any $g\in\Gamma_{s}$, we have $(g(B^{+}_{2})\cap\partial B_{\rho})\backslash W_{\varepsilon_{1}}\neq \emptyset$, using $(f_{1})$, we have
\begin{align*}
\inf_{g\in\Gamma_{s}}\sup_{g(B_{2}^{+})\backslash W_{\varepsilon_{2}}}I(u)&\geq\inf_{\partial B_{\rho}}I \\
&=\inf_{\|u\|=\rho}\frac{1}{2}\|u\|^{2}-\int F(u) \\
&\geq\inf_{\|u\|=\rho}\frac{1}{4}\|u\|^{2}-C\|u\|^{p}>\alpha>0 \\
\end{align*}
for some suitable positive numbers $\alpha$ and $\rho$. Suppose $K_{c_{s}}\backslash W_{\varepsilon_{2}}=\emptyset$, we have a continuous mapping $\eta:[0,1]\times E\to E$ such that
\begin{itemize}
  \item[$(1)$] $\eta(0,u)=u$ for any $u\in E$.
  \item[$(2)$] $\eta(t,u)=u$ for all $t\in[0,1]$,$\,\,u\in (I_{c-2\varepsilon}^{c+2\varepsilon})^{c}\cup W_{\varepsilon_{1}}$.
  \item[$(3)$] $\eta(1,I^{c_{s}+\varepsilon}\backslash W_{\varepsilon_{2}})\subset I^{c_{s}-\varepsilon}$.
  \item[$(4)$] $\eta(t,\overline{P_{\varepsilon_{2}}^{\pm}})\subset\overline{P_{\varepsilon_{2}}^{\pm}}$ for all $t\in[0,1]$
\end{itemize}
due to Lemma \ref{l:mapping2}. Let $\varepsilon=\frac{c_{s}-\alpha}{4}>0$. Select a $g_{0}\in\Gamma_{s}$ satisfies
$$\sup_{g_{0}(B^{+}_{2})\backslash W_{\varepsilon_{2}}} I<c_{s}+\varepsilon.$$
Denote $g_{1}=\eta\circ g_{0}$. Then we have
\begin{align*}
g_{1}(B_{2}^{+})\backslash W_{\varepsilon_{2}} & = \eta\circ g_{0}(B_{2}^{+})\backslash W_{\varepsilon_{2}} \\
&=\eta\big((g_{0}(B_{2}^{+})\backslash W_{\varepsilon_{2}})\cup W_{\varepsilon_{2}}\big)\backslash W_{\varepsilon_{2}}\\
&\subset\eta\big(g(B_{2}^{+})\backslash W_{\varepsilon_{2}}\big)\subset \eta(I^{c_{s}+\varepsilon})\subset I^{c_{s}-\varepsilon}\\
\end{align*}
and
\begin{itemize}
  \item $g_{1}|_{\partial_{0}B_{2}^{+}}=\eta\circ g_{0}|_{\partial_{0}B_{2}^{+}}=id$ since $\eta=id$ in $\big(I_{c-2\varepsilon}^{c+2\varepsilon}\big)^{c}$;
  \item $g_{1}|_{B_{2}^{+}\cap \mathbb{R}e_{1}}=\eta\circ g_{0}|_{B_{2}^{+}\cap \mathbb{R}e_{1}}=id$ since $\eta=id$ in $W_{\varepsilon_{1}}$ and $B_{2}^{+}\cap \mathbb{R}e_{1} \subset W_{\varepsilon_{1}}$.
\end{itemize}
Therefore we can find a $g_{1}\in\Gamma_{s}$ with
$$\sup_{g_{1}(B_{2}^{+})\backslash W_{\varepsilon_{2}}}I\leq c_{s}-\varepsilon.$$
Then, we have a contradiction with the definition of $c_{s}$.
\begin{flushright}
\qed
\end{flushright}

\section{Remarks}
\subsection{More constructions of the sign-changing minimax values}

In this section, we will give two more constructions of the sign-changing minimax values, and the second one is a slight modification of the one in \cite{LLW}. Since we will not use Lemma \ref{l:mapping2} in this section, which means we do not use different neighbourhoods of the positive and negative cones, we denote $W_{\varepsilon}$ by $W$ for short. And always assume that the parameter $\varepsilon$ is small enough. Define the following two minimax values. The first one is:
$$c_{s}':=\inf_{g\in\Gamma_{s}'}\sup_{g(B_{2}^{+})\backslash W}I,$$
where
\begin{itemize}
  \item $\Gamma_{s}'=\big\{g\in C(B_{2}^{+},E)\big|\,\,g|_{\partial_{0}B_{2}^{+}}=id,\,\,g(\partial_{1}B_{2}^{+})\subset \overline{P_{\varepsilon}^{+}}\,\,and\,\,g(\partial_{2}B_{2}^{+})\subset \overline{P_{\varepsilon}^{-}}\big\}$;
  \item $B_{2}^{+}=\big\{xe_{1}+ye_{2}|\,\,y\geq0\,\,and\,\,x^{2}+y^{2}\leq R^{2}\big\}$;
  \item $\partial_{0}B_{2}^{+}=\{x^{2}+y^{2}=R^{2}\}\cap B_{2}^{+}$, $\partial_{1}B_{2}^{+}=\{x\geq0\,\,and\,\,y=0\}\cap B_{2}^{+}$ and $\partial_{2}B_{2}^{+}=\{x\leq0\,\,and\,\,y=0\}\cap B_{2}^{+}$,
\end{itemize}
and the functions $e_{1}$, $e_{2}$ and the number $R$ are the same with Section 3. The second one is:
$$c_{s}'':=\inf_{g\in\Gamma_{s}''}\sup_{g(D^{++}_{2})\backslash W}I,$$
where
\begin{itemize}
  \item $\Gamma_{s}''=\big\{g\in C(D_{2}^{++},E)\big|\,\,g|_{\partial_{0}D_{2}^{++}}=id,\,\,g(\partial_{1}D_{2}^{++})\subset \overline{P_{\varepsilon}^{+}}\,\,and\,\,g(\partial_{2}D_{2}^{++})\subset \overline{P_{\varepsilon}^{-}}\big\}$;
  \item $D_{2}^{++}=\{x\alpha_{1}+y\alpha_{2}|\,\,x,y\geq 0 \,\,and\,\, x^{2}+y^{2}\leq R^{2}\}$;
  \item $\partial_{0}D_{2}^{++}=\{x^{2}+y^{2}=R^{2}\}\cap D_{2}^{++}$, $\partial_{1}D_{2}^{++}=\{x\geq0\,\,and\,\,y=0\}\cap D_{2}^{++}$ and $\partial_{2}D_{2}^{++}=\{x=0\,\,and\,\,y\geq0\}\cap D_{2}^{++}$;
  \item $\alpha_{1},\,\,\alpha_{2}\in C_{0}^{\infty}(\W)$ are two nonzero functions with unit Sobolev norm and $\alpha_{1}\geq 0$, $\alpha_{2}\leq 0$ on $\W$ and $supp(\alpha_{1})\cap supp(\alpha_{2})=\emptyset$.
\end{itemize}
In the following paragraph, it only need to check the linking results. The proof of the existence of sign-changing critical point is the same as we did in Section 3. We only verify the linking result for the first minimax value since the second one is similar.
\begin{proposition}
For any $g\in\Gamma_{s}'$,
$$\big(g(B_{2}^{+})\cap\partial B_{\rho}\big)\backslash W\neq\emptyset.$$
\end{proposition}
\noindent{\bf Proof.}
The proof is divided into the following steps.
\subparagraph{Step 1. Modification near the origin.}
In the whole process, we always assume that $\varepsilon>0$ is small enough. Since the mapping $g\in\Gamma_{s}'$, $g(0)\in\overline{P_{\varepsilon}^{+}}\cup\overline{P_{\varepsilon}^{-}}$, then there are $\rho_{1},\,\,\rho_{2}>0$ such that
$$g(B_{2}^{+}\cap B_{\rho_{1}})\subset B_{\rho_{2}}$$
and $\rho_{2}$ are small enough, say $10\rho_{2}<\rho$. Then we can modify $g$ on the set $B_{2}^{+}\cap B_{\rho_{1}}$ by
\begin{equation}
\overline{g}(u)=\left\{
\begin{aligned}
&g(u)\qquad  u\in B_{2}^{+}\backslash B_{\rho_{1}};\nonumber \\
&\frac{\|u\|}{\rho_{1}}g\big(\rho_{1}\frac{u}{\|u\|}\big) \qquad u\in \big(B_{2}^{+}\cap B_{\rho_{1}}\big)\backslash \{0\};\nonumber \\
&0\qquad   u=0.\nonumber
\end{aligned}
\right.
\end{equation}
It is easy to verify that $\overline{g}(0)=0$,
$$\overline{g}(B_{2}^{+})\cap \partial B_{\rho}=g(B_{2}^{+})\cap\partial B_{\rho},$$
$$\overline{g}(\partial_{1}B_{2}^{+})\subset\overline{P_{\varepsilon}^{+}}$$
and
$$\overline{g}(\partial_{2}B_{2}^{+})\subset\overline{P_{\varepsilon}^{-}},$$
since $\overline{P_{\varepsilon}^{+}}$ and $\overline{P_{\varepsilon}^{-}}$ are convex.

\subparagraph{Step 2. Modification in $\overline{P_{\varepsilon}^{-}}$.}
To use the genus argument, we need to modify $\overline{g}$ by a mapping whose restriction on $\partial_{1}B_{2}^{+}\cup\partial_{2}B_{2}^{+}$ is odd. To do this, we define a homeomorphism
$$\phi:B_{2}^{+}\to\overline{B_{2}^{+}\backslash C}=:D,$$
where
$$C=\big\{-R\leq x\leq 0\,\,and\,\,0\leq y\leq \frac{R}{3}|sin\frac{\pi x}{R}|\big\}\cap B_{2}^{+},$$
and
\begin{itemize}
  \item $\phi|_{\partial_{0}B_{2}^{+}}=id$;
  \item $\phi|_{\partial_{1}B_{2}^{+}}=id$;
  \item On $\partial_{2}B_{2}^{+}$, $\phi(xe_{1})=xe_{1}+\frac{R}{3}|sin\frac{\pi x}{R}|$ and hence $\phi(\partial_{2}B_{2}^{+})=D\cap C$.
\end{itemize}
Then we define the following mapping $\overline{\overline{g}}$. For the sake of
convenient, we express it into coordinate form.
\begin{equation}
\overline{\overline{g}}(xe_{1}+ye_{2})=\left\{
\begin{aligned}
&\overline{g}\circ\phi^{-1}(xe_{1}+ye_{2}), \qquad x,y\in D;\nonumber \\
&-\overline{g}(-xe_{1}), \qquad u\in B_{2}^{+}\cap \{x\leq 0 \,\,and\,\,y=0\};\nonumber \\
&-\Big(1-\frac{3y}{R|sin\frac{\pi x}{R}|}\Big)\overline{g}(-xe_{1})+\qquad \nonumber\\
&\frac{3y}{R|sin\frac{\pi x}{R}|}\overline{g}\circ\phi^{-1}(xe_{1}+\frac{R}{3}|sin\frac{\pi x}{R}|e_{2}),\quad   u\in C\backslash \big(D\cup\{y=0\}\big).\nonumber
\end{aligned}
\right.
\end{equation}
It is easy to verify that
\begin{itemize}
  \item $\overline{\overline{g}}$ is continuous on $B^{+}_{2}$;
  \item $\overline{\overline{g}}\big|_{\partial_{0}B_{2}^{+}}=id$;
  \item $\overline{\overline{g}}$ is odd on $B_{2}^{+}\cap\{y=0\}=\partial_{1}B_{2}^{+}\cup\partial_{2}B_{2}^{+}$.
\end{itemize}

\subparagraph{Step 3. Odd extension.}
Define
\begin{equation}
\overline{\overline{\overline{g}}}(u)=\left\{
\begin{aligned}
\overline{\overline{g}}(u),\qquad  & u\in B_{2}^{+}; \nonumber\\
-\overline{\overline{g}}(-u),\qquad & -u\in B_{2}^{+}. \nonumber
\end{aligned}
\right.
\end{equation}
Using the genus argument, we have
$$\big(\overline{\overline{\overline{g}}}(B_{2}^{+})\cap \partial B_{\rho}\big)\backslash W\neq\emptyset,$$
which implies that
$$\big(\overline{\overline{g}}(B_{2}^{+})\cap \partial B_{\rho}\big)\backslash W\neq\emptyset.$$
On one hand, since
\begin{itemize}
  \item $\overline{\overline{g}}(C)\subset\overline{P_{\varepsilon}^{-}}$ since $\overline{\overline{g}}(C\cap D)\subset\overline{P_{\varepsilon}^{-}}$, $\overline{\overline{g}}(C\cap\{y=0\})\subset\overline{P_{\varepsilon}^{-}}$ and $\overline{P_{\varepsilon}^{-}}$ is convex;
  \item $\phi$ is a homeomorphism,
\end{itemize}
we have
\begin{align*}
\overline{\overline{g}}(B_{2}^{+})\backslash W & = \overline{\overline{g}}(D\cup C)\backslash W\\
&=\big(\overline{\overline{g}}(D)\cup\overline{\overline{g}}(C)\big)\backslash W =\overline{\overline{g}}(D)\backslash W\\
&=\overline{g}(\phi^{-1}(D))\backslash W=\overline{g}(B_{2}^{+})\backslash W.
\end{align*}
On the other hand, recall the result in Step 1,
$$\overline{g}(B_{2}^{+})\cap \partial B_{\rho}=g(B_{2}^{+})\cap\partial B_{\rho},$$
we have
\begin{align*}
\big(g(B_{2}^{+})\cap\partial B_{\rho}\big)\backslash W & = \big(g(B_{2}^{+})\backslash W \big)\cap\partial B_{\rho}\\
&=\big(\overline{g}(B_{2}^{+})\backslash W \big)\cap\partial B_{\rho}\\
&=\big(\overline{\overline{g}}(B_{2}^{+})\backslash W \big)\cap\partial B_{\rho}\\
&=\big(\overline{\overline{g}}(B_{2}^{+})\cap\partial B_{\rho}\big)\backslash W\neq\emptyset,
\end{align*}
which is what we want.
\begin{flushright}
\qed
\end{flushright}

Following a similar process, we can conclude the linking result for the second minimax value.
\begin{proposition}
For any $g\in\Gamma_{s}''$, we have
$$(g(D_{2}^{++})\cap\partial B_{\rho})\backslash W\neq \emptyset.$$
\end{proposition}
Then we can claim that $K_{c_{s}'}\backslash W\neq\emptyset$ and $K_{c_{s}''}\backslash W\neq \emptyset$.

\subsection{A remark on the linking sets}
In the proceeding paragraph, when we deal with the linking structure of sign-changing critical points, we used the 1-codimensional set $\partial B_{\rho}$. In fact, we can use some sets of higher codimension. For example, consider the set $F=Y\cap\partial B_{\rho}$, where the subspace $Y$ satisfies that $codim(Y)=1$ and $F\cap(E\backslash W)=\emptyset$. It is valid for all three kinds of critical values, we only consider the first kind as an example. Using the intersection lemma (c.f. \cite[Lemma 6.4]{S1}), we have for any $g\in\Gamma_{s}$,
$$g(B_{2}^{+})\cap F\neq\emptyset,$$
and hence
$$\big(g(B_{2}^{+})\cap F\big)\backslash W\neq\emptyset.$$
In this case, $codim(F)=2$. It should be remarked that in \cite{LLW}, the linking set is $\partial\overline{P_{\varepsilon}^{+}}\cap\partial\overline{P_{\varepsilon}^{-}}$. These examples show that the set $\partial B_{\rho}$ is not the 'smallest' linking set. It could be intersecting to look for a necessary and sufficient condition for a set to be a linking set for sign-changing mountain pass lemmas.

\section*{Acknowledgement}
The author would like to express his sincere thanks to his advisor Professor Zhi-Qiang Wang for useful suggestions, infinitely patience and language help. This research is supported by CNSF(11771324).

%%%%%%%%%%%%%%%%%%%%%%%%%%%%%%%%%%%%%%%%%%%%%%%%%%%%%%%%%

\end{document}